\documentclass[12pt]{article}
\usepackage{amsfonts,latexsym}

\newcommand{\PP}{{\bf P}}

\newcommand{\Proof}{{\it Proof\/.\quad }}
\newcommand{\F}{{\rm I\!F}}

\newcommand{\D}{{\rm I\!D}}
\newcommand{\R}{{\rm I\!R}}

\newcommand{\N}{{\rm I\!N}}

\title{ On solutions of equations with measurable coefficients driven by $\alpha$- stable processes}
\author{V.~P.~KURENOK$$\\
	$$Department of Electrical and Systems Engineering,\\
	Washington University in St. Louis,\\
	One Brookings Drive, St. Louis, MO 63130-4899, USA\\
	e-mail: kurenokv@ese.wustl.edu}

\begin{document}

\newtheorem{Th}{Theorem}[section]
\newtheorem{Def}[Th]{Definition}
\newtheorem{Prop}[Th]{Proposition}
\newtheorem{Le}[Th]{Lemma}
\newtheorem{Cor}[Th]{Corollary}
\newtheorem{Rems}[Th]{Remarks}
\newtheorem{Rem}[Th]{Remark}
\newtheorem{Ex}[Th]{Example}
\thispagestyle{empty}
\maketitle
\vspace*{0.5cm}

\begin{abstract}

We prove the existence of solutions for the stochastic differential equation  $dX_t=b(t,X_{t-})dZ_t+a(t,X_t)dt, X_0\in\R, t\ge 0,$  with only measurable coefficients $a$ and $b$ satisfying the condition $0<\mu\le |b(t,x)|\le \nu$ and $|a(t,x)|\le K$ for all $t\ge 0, x\in\R$ where $\mu, \nu, $ and $K$ are some constants. The driving process $Z$ is a symmetric stable process of index $1<\alpha<2$. This generalizes the result of N. V. Krylov \cite{Krylov} for the case of $\alpha=2$, that is when $Z$ is a Brownian motion. The proof is based on integral estimates of Krylov type for the given equation which are also derived in the note and are of independent interest. Moreover, unlike in \cite{Krylov}, we use a different approach to derive the corresponding integral estimates.   
\end{abstract}
{\small{\it AMS Mathematics subject classification. Primary}\hspace{15pt} 60H10,
60J60, 60J65, 60G44}

{\small{\it Keywords and phrases.}\hspace{15pt} Stochastic differential equations, symmetric stable processes, Krylov's estimates, Fourier transform}

\section{Introduction}

\quad We consider here a stochastic differential equation of the form
\begin{equation}\label{eq1}
dX_t=b(t,X_{t-})dZ_t+a(t,X_t)dt, X_0=x_0\in\R, t\ge 0.
\end{equation}
The existence of solutions for equation (\ref{eq1}) with {\it only measurable coefficients} $a$ and $b$ was proved first by N. V. Krylov in \cite{Krylov} for the case when the driving process $Z$ is a Brownian motion. The proof was based on using of corresponding integral estimates for solutions $X$ of (\ref{eq1}) he was also first to derive. These integral estimates turned later to be very useful  in various areas of stochastic processes including the optimal control of processes described by equation (\ref{eq1}). The estimates of such kind are now often refereed to as {\it Krylov's estimates}.

\quad In order to prove the corresponding integral estimates, Krylov had used the Bellman principle of optimality known in the control theory of stochastic processes. For that, given a smooth function $f(t,x)$, he considered a  value function
\begin{equation}\label{valuefunction}
v(t,x):=\sup_{\beta\in{\cal B}}{\bf E}\int_0^{\infty}e^{-\phi_s^{\beta}}\psi_s^{\beta}f(t+r^{\beta}_s,x+X_s^{\beta})ds
\end{equation}
where $(\phi^{\beta},\psi^{\beta})$  and $(r^{\beta},X^{\beta})$ are appropriately chosen stochastic processes and ${\cal B}$ is a suitably chosen set of control parameters.
	
\quad Using (\ref{valuefunction}), one derived then the corresponding Bellman equation for the function $v(t,x)$ and upon integrating it one received the estimates of the form 
\begin{equation}\label{supnormestimate}
\sup_{t,x} v(t,x)\le M\|f\|_{L_p},
\end{equation}
where $\|f\|_{L_p}$ is the $L_p$-norm of the function $f$ and $p\ge 1$. Finally, using the It{\'o}'s formula and the estimates (\ref{supnormestimate}), one obtained the integral  estimates
\begin{equation}\label{krylovestimate}
{\bf E}\int_0^{\infty}f(s,X_s)ds\le M\|f\|_{L_p},
\end{equation}
known now as Krylov's estimates.

\quad As an application of (\ref{krylovestimate}), Krylov proved the existence of solutions of equation (\ref{eq1}) in the case when $Z$ is a Brownian motion and the measurable coefficients $a(t,x)$ and $b(t,x)$ are such that, for all $(t,x)$, it holds
\begin{equation}\label{eq2}
 0<\mu\le |b(t,x)|\le \nu, \mbox{ } |a(t,x)|\le K
\end{equation}
for some constants $\mu,\nu,$ and $K$.

\quad In this note we consider the equation (\ref{eq1}) when the driving process $Z$ is a symmetric stable process of index $1<\alpha\le 2$. For $\alpha=2$, $Z$ is then a Brownian motion process.

\quad  One of the main results  here is the proof of the existence of solutions of equation (\ref{eq1}) when the coefficients $a$ and $b$ are only measurable and satisfy the condition (\ref{eq2}). This extends the result of Krylov for the Brownion motion case to case of a symmetric stable process with the index $1<\alpha\le 2$.

\quad To prove the existence of solutions, we will first derive the corresponding Krylov's estimates for processes $X$ satisfying (\ref{eq1}). However, in order to do so, unlike in \cite{Krylov}, we do not use any facts from the optimal control theory for stochastic processes but consider a parabolic integro-differential equation of the form 
\begin{equation}\label{eq3}
u_t+|b|^{\alpha}{\mathcal L}u +au_x-\lambda(1+|b|^{\alpha})u+f=0
\end{equation}
where ${\cal L}$ is the generator of the process $Z$ (see definitions below) and $\lambda>0$. 

\quad Assuming that the functions $a$ and $b$ satisfy the condition (\ref{eq2}), we will prove some important {\it a priori estimates} for the equation (\ref{eq3}) of the form
\begin{equation}\label{Hnormestimate}
\|u\|_H\le M\|f\|_{L_2}
\end{equation}
which, in turn, will imply the estimates 
\[
\sup_{t,x} u(t,x)\le M\|f\|_{L_2}.
\]
\quad Moreover, a priori estimates (\ref{Hnormestimate}) are then also used to prove the existence of a solution $u$ of the equation (\ref{eq3}) given a fixed function $f\in L_2$. The later fact is important to derive the Krylov's estimates.

\begin{Rem}\label{BE} Following Krylov's idea in the Brownian motion case, one of the possible Bellman equations related to equation (\ref{eq1}) can be derived as follows. One considers the controlled process $(t+r_s^{\beta},x+X_s^{\beta})$ defined as 
\[
dr_s^{\beta}=(1-|\sigma_s|^{\alpha})ds, r_0=0,
\]
\[
dX_s^{\beta}=\gamma_sds+\sigma_sdZ_s, X_0=0,
\]
with the value function
\[
v(t,x)=\sup_{{\cal B}}{\bf E}\int_0^{\infty}e^{-\lambda s}\psi_s^{\beta}f(t+r_s^{\beta},x+X_s^{\beta})ds,
\]
where $f\in C_0^{\infty}(\R^2)$\footnote{$C_0^{\infty}(\R^2)$ defines the class of infinitely differentiable functions with a compact support on $\R^2$} and $\lambda>0$.

Here ${\cal B}$ is the class of strategies $\beta_s=(\gamma_s,\sigma_s)$ such that
\[
|\sigma_s|\le 1, |\gamma_s|\le K|\sigma_s|^{\alpha}
\]
and
\[
\psi_s^{\beta}=\sqrt{(1-|\sigma_s^{\alpha}|)|\sigma_s|^{\alpha}}.
\]

The corresponding Bellman equation ($\sigma$ and $\gamma$ are numbers) holds a.e. in $\R^2$:
\[
\sup_{|\sigma|\le 1}\sup_{|\gamma|\le K|\sigma|^{\alpha}}\Bigl[(1-|\sigma|^{\alpha})v_t+|\sigma|^{\alpha}{\cal L}v-\lambda v+\gamma v_x+\sqrt{(1-|\sigma|^{\alpha})|\sigma|^{\alpha}}f\Bigr]=0,
\]
where $v_t$ and $v_x$ are partial derivatives of $v$ in $t$ and $x$, respectively.

It is then not hard to see that the Bellman equation is equivalent to
\begin{equation}\label{Bellman}
(v_t-\lambda v)({\cal L}v-\lambda v+K|v_x|)=\frac{1}{4}f^2
\end{equation}
which also holds a.e. in $\R^2$. However, it has been unclear to us how to integrate  the equation (\ref{Bellman}) to prove the estimates (\ref{Hnormestimate}). This lead us to consider an alternative route in form of the equation (\ref{eq3}).
\end{Rem}

\quad Finally, we give a brief overview of existence results for the equation (\ref{eq1}) with only measurable coefficients $a$ and $b$ and $1<\alpha<2$ known for some particular cases.

\quad The equation (\ref{eq1}) without drift (that is, when $a=0$), and time-independent coefficient  $b(x)$ was studied in detail by P. A. Zanzotto  \cite{Za1} where he relied on the systematic use of time change techniques.

\quad The time-dependent equation (\ref{eq1}) without drift was studied by H. Pragarauskas and P. A. Zanzotto \cite{PragZa}. To prove the existence of solutions, one used the method of integral estimates similar to \cite {Krylov}. The corresponding integral estimates were proven by H. Pragarauskas in \cite{Prag}. Some other sufficient existence conditions for the time-dependent equation without drift different from those in \cite{PragZa} and with $0<\alpha<2$ were found in \cite{EngKu}. 

\quad The equation (\ref{eq1}) with time-independent coefficients $a(x)$ and $b(x)$  was considered by the author in \cite{Ku1} where he proved the existence of solutions with only measurable coefficients $a$ and $b$ satisfying the condition
\[
0<\mu\le|b(x)|\le\nu, \mbox { }, |a(x)|\le K|b(x)|^{\alpha}.
\] 
for all $x\in\R$ and some constants $\mu, \nu,$ and $K$.

\section{Some preliminary facts}
\setcounter{equation}{0}

\quad By ${\bf D}_{[0,\infty)}(\R)$ we denote, as usual, the Skorokhod space, i.e. the set of all real-valued functions $z:[0,\infty)\to \R$ with right-continuous trajectories and with finite left limits (also called {\it c\'{a}dl\'{a}g} functions). For simplicity, we shall write ${\bf D}$ instead of ${\bf D}_{[0,\infty)}(\R)$. We will equip ${\bf D}$ with the $\sigma$-algebra $\mathcal{D}$ generated by the Skorokhod topology. Under ${\bf D}^n$ we will understand the $n$-dimensional Skorokhod space defined as ${\bf D}^n={\bf D}\times\dots\times{\bf D}$ with the corresponding $\sigma$-algebra $\mathcal{D}^n$ being the direct product of $n$ one-dimensional $\sigma$-algebras $\mathcal{D}$.  

\quad Let $(\Omega, \mathcal{F},\PP)$ be a complete probability space with a filtration $\F=(\mathcal{F}_t)$. We use the notation $(Z,\F)$ to indicate that a process  $Z$ is adapted to $\F$. A process $(Z, \F)$ is called a symmetric stable process of index $\alpha\in (0,2]$ if trajectories of $Z$ are
c\'{a}dlag functions and ${\bf
E}\left(\exp\left(i\xi(Z_t-Z_s)\right)|\mathcal{F}_s\right)=\exp\left(-(t-s)c|\xi|^{\alpha}\right)$
for all $t>s\ge 0$ and $\xi\in{\R}$, where $c>0$ is a constant. The function
$\psi(\xi)=c|\xi|^{\alpha}$ is called the characteristic exponent of the
process $Z$.

\quad The process $Z$ is a process with independent increments thus a
Markov process. Therefore, it can be characterized in terms of Markov processes.
For any function $g\in L^{\infty}(\R)$ and $t\ge 0$, define the operator
\[
(T_tg)(x):={\bf E}\Bigl(g(x+Z_t)\Bigr)
\]
where $L^{\infty}(\R)$ is the Banach space of functions $g:\R\to\R$ with the
norm $\|g\|_{\infty}=ess\sup|g(x)|$. For a suitable class of functions $g(x)$, we can define an operator $\mathcal{L}$ often called the infinitesimal generator of the process $Z$ as
\begin{equation}\label{infgenerator}
(\mathcal{L}g )(x)=\lim_{t\downarrow 0}\frac{(T_tg)(x)-g(x)}{t}.
\end{equation}
\quad On another hand, in the case of $\alpha\in(0,2)$, the process $Z$ as a
purely discontinuous Markov process can be described by its Poisson jump measure
(jump measure of $Z$ on interval $[0,t]$) defined as
$$
N(U\times[0,t])=\sum_{s\le t}1_{U}(Z_s-Z_{s-}),
$$
the number of times before the time $t$ that $Z$ has jumps whose size lies in
the set $U$. The compensating measure of $N$, say $\hat N$, is given by
\[
\hat N(U)={\bf E}N(U\times[0,1])=\int_{U}\frac{1}{|x|^{1+\alpha}}dx.
\]
It is known that for $\alpha<2$ 
\begin{equation}\label{generatorform}
(\mathcal{L}g)(x)=\int\limits_{\R\setminus\{0\}}[g(x+z)-g(x)-{\bf
1}_{\{|z|<1\}}g^{\prime}(x)z]\frac{c_1}{|z|^{1+\alpha}}dz
\end{equation}
for any $g\in C^2_b(\R)$, where $C^2_b(\R)$ is the set of all bounded and twice continuously
differentiable functions $g:\R\to\R$ whose derivatives are also bounded. We shall assume from now on the constant
$c_1$ to be chosen in the way that $\psi(\xi)=1/2|\xi|^{\alpha}$. In the case of $\alpha=2$ the infinitesimal generator of $Z$ is the second derivative operator, that is, ${\mathcal
L}g(x)=\frac{1}{2}g^{\prime\prime}(x)$.

\quad Let $L_p(\R^2), p\ge 1,$ define the space of all measurable functions
$g:\R^2\to\R$ such that $(\int_{\R^2}|g(s,x)|^pdsdx)^{1/p}<\infty$.
Then, for any $g\in L_1(\R^2)$, there exists its Fourier transform $Fg$ defined
as  
\[
[Fg](\tau,\omega):=\int\limits_{\R^2}
e^{is\tau}e^{ix\omega}g(s,x)dsdx,\quad (\tau,\omega)\in\R^2.
\]
Moreover,  if $Fg\in L_1(\R^2)$, then also the inverse Fourier transform
$F^{-1}$ of the function $Fg$ exists and
\begin{equation}\label{Fourierinverse}
g(s,x)=\frac{1}{(2\pi)^2}\int_{\R^2}[Fg](\tau,\omega)e^{-is\tau}e^{-ix\omega}d\tau d\omega,\quad
(s,x)\in\R^2.
\end{equation}
\quad We also note that calculating the Fourier transform of a function of two
variables can be performed as calculating the single Fourier transform in one
variable and then in another, in any order. The next statement is known (see, for example, \cite{Be},  Proposition 9, ch. 1) but we provide a short proof of it for the convenience of the reader. 

\begin{Prop}\label{Fouriertransform} Assume $0<\alpha\le 2$. The following
statements are true:
\begin{enumerate}
\item[(i)] For every function $g\in L^{\infty}\cap L_1(\R)$, it holds
\[
F(T_tg)=e^{-\frac{1}{2}t|\omega|^{\alpha}}Fg.
\]
\item[(ii)] Assume $g\in C^{\infty}_0(\R)$ and $\mathcal{L}g\in L_1(\R)$. Then
\[
F(\mathcal{L}g)=-\frac{1}{2}|\omega|^{\alpha}Fg.
\]

\end{enumerate}
\end{Prop}
\Proof For $(i)$:
\[
F(T_tg)(\omega)=\int\limits_{\R}e^{i\omega x}T_tg(x)dx={\bf
E}\Bigl(\int\limits_{\R}e^{i\omega x}g(x+Z_t)dx\Bigr)=
\]
\[
{\bf E}\Bigl(\int\limits_{\R}e^{i\omega (y-Z_t)}g(y)dy\Bigr)={\bf E}e^{-i\omega
Z_t}\int\limits_{\R}e^{i\omega y}g(y)dy=e^{-t|\omega|^{\alpha}}Fg(\omega).
\]
The statement $(ii)$ follows from $(i)$ and the definition (\ref{infgenerator}):
\[
F(\mathcal{L}g)(\omega)=\int\limits_{\R}e^{i\omega
x}\mathcal{L}g(x)dx=\int\limits_{\R}e^{i\omega x}\lim_{t\downarrow
0}\frac{T_tg(x)-g(x)}{t}dx=
\]
\[
\lim_{t\downarrow 0}\frac{1}{t}\Bigl(\int\limits_{\R}e^{i\omega
x}T_tg(x)dx-\int\limits_{\R}e^{i\omega x}g(x)dx\Bigr)=
\]
\[
\lim_{t\downarrow
0}\frac{1}{t}\Bigl(e^{-t\frac{1}{2}|\omega|^{\alpha}}Fg(\omega)-Fg(\omega)\Bigr)=
\]
\[
Fg(\omega)\lim_{t\downarrow
0}\frac{e^{-t\frac{1}{2}|\omega|^{\alpha}}-1}{t}=-\frac{1}{2}|\omega|^{\alpha}Fg(\omega).
\]
$\Box$

\quad We also introduce the following space of functions associated with the infinitisimal  operator ${\cal L}$ of a symmetric stable process of index $\alpha$. For any $u\in C_0^{\infty}(\R^2)$, define the norm
\begin{equation}\label{sobolevnorm}
\|u\|_H:=\|u\|_{L_2}+\|u_t\|_{L_2}+\|{\mathcal L}u\|_{L_2}.
\end{equation}

We say that a function $u(t,x)\in L_2(\R^2)$ belongs to the space $H(\R^2)$ if there is a sequence of functions $u^n\in C_0^{\infty}(\R^2)$ such that $\|u^n\|_H<\infty$ for all $n=1,2,..$,
\[
\|u^n-u\|_{L_2}\to 0
\]
as $n\to\infty$, and
 \[
\|u_t^n-u_t^m\|_{L_2}\to 0, \|u^n-u^m\|_{L_2}\to 0, \|{\cal L}u^n-{\cal L}u^m\|_{L_2}\to 0
\]
as $n,m\to\infty$. The space $H$ is then called a {\it Sobolev space}.

\section{ Analytic a priori estimates} 
\setcounter{equation}{0}

\quad Let $\lambda>0$ and $\alpha\in(1,2)$. In this section we consider the integro-differential equation of parabolic type (\ref{eq3}) in the Sobolev space $H$ 
with the norm $\|\cdot\|_H$ defined in (\ref{sobolevnorm}). Moreover, we choose a function $f\in C_0^{\infty}(\R^2)$ and assume that the coefficients $b(t,x)$ and $a(t,x)$ satisfy the condition (\ref{eq2}).

\quad We are interested in deriving some {\it a priori estimates} for a solution $u$ of the equation (\ref{eq3}) in terms of the $L_2$-norm of the function $f$. Since the existence of a solution is not known yet, such estimates are called a priori estimates. These estimates are crucial for deriving the integral estimates of Krylov type for processes $X$ satisfying the stochastic equation (\ref{eq1}).  

\quad Moreover, those a priori estimates derived here can be used to actually prove the existence of a solution $u\in H$ of equation (\ref{eq3}) for any $f\in L_2$. The corresponding proof is based on the {\it method of continuity} and the {\it method of a priori estimates} known in the theory of classical elliptic and parabolic equations, that is when ${\cal L}$ is the second derivative operator. The proof of the existence of a solution of equation (\ref{eq3}) is provided in the Appendix.

\begin{Le}\label{estimate1} Let $u\in C^{\infty}_0(\R^2)$ be a solution of the equation (\ref{eq3}). Then, it holds
\begin{equation}\label{eq4}
\|u\|_H\le M\|f\|_{L_2}.
\end{equation}
\end{Le}
\Proof It follows from (\ref{eq3}) that
$$
[(u_t-\lambda u)+|b|^{\alpha}({\mathcal L}u-\lambda u)]^2=(au_x+f)^2\le
2a^2u_x^2+2f^2
$$
and
$$
\frac{1}{|b|^{\alpha}}(u_t-\lambda u)^2+2(u_t-\lambda u)({\mathcal L}u-\lambda
u)+|b|^{\alpha}({\mathcal L}u-\lambda
u)^2\le\frac{2}{|b|^{\alpha}}(K^2u_x^2+f^2).
$$
The condition (\ref{eq2}) implies that
\begin{equation}\label{eq5}
\frac{1}{\nu^{\alpha}}(u_t-\lambda u)^2+2(u_t-\lambda u)({\mathcal L}u-\lambda
u)+\mu^{\alpha}({\mathcal L}u-\lambda u)^2\le\frac{2}{\mu^{\alpha}}(K^2u_x^2+f^2).
\end{equation}

Using the Plansherel's identity and Proposition \ref{Fouriertransform}, we obtain
\begin{equation}\label{eq6}
\int_{\R^2}(u_t-\lambda u)^2dtdx=\int_{\R^2}|F(u_t-\lambda
u)|^2d\tau dw=\int_{\R^2}|F(u)|^2(\lambda^2+\tau^2)d\tau d\omega,
\end{equation}
\begin{equation}\label{eq7}
\int_{\R^2}({\mathcal L}u-\lambda u)^2dtdx=\int_{\R^2}|F({\mathcal
L}u-\lambda
u)|^2d\tau dw=\int_{\R^2}|F(u)|^2(\lambda+|\omega|^{\alpha})^2d\tau d\omega,
\end{equation}
and
\begin{equation}\label{eq8}
\int_{\R^2}u_x^2dtdx=\int_{\R^2}|\omega|^2|F(u)|^2d\tau d\omega.
\end{equation}
As it can be easily seen, there is $\delta>0$ so that
\begin{equation}
\label{eq9}
\mu^{\alpha}(\lambda+|\omega|^{\alpha})^2\ge \frac{4K^2}{\mu^{\alpha}}|\omega|^2
\end{equation}

for all $\omega\in\R$ and all $\lambda\ge \delta$.

Now, we integrate equation (\ref{eq5}) over $\R^2$ and use identities
(\ref{eq6})-(\ref{eq9}) to obtain
$$
\frac{1}{\nu^{\alpha}}\int_{\R^2}|F(u)|^2(\lambda^2+\tau^2)+2\int_{\R^2}(u_t-\lambda
u)({\mathcal L}u-\lambda
u)+\frac{\mu^{\alpha}}{2}\int_{\R^2}(\lambda+|\omega|^{\alpha})^2|F(u)|^2
$$
\begin{equation}\label{eq50}
\le \frac{2}{\mu^{\alpha}}\int_{\R^2}f^2.
\end{equation}
The last inequality implies
$$
\frac{\lambda^2}{\nu^{\alpha}}\int_{\R^2}|F(u)|^2+2\int_{\R^2}(u_t-\lambda u)({\mathcal
L}u-\lambda u)+\frac{\mu^{\alpha}\lambda^2}{2}\int_{\R^2}|F(u)|^2\le
\frac{2}{\mu^{\alpha}}\int_{\R^2}f^2,
$$
or
\begin{equation}\label{eq10}
(\frac{\mu^{\alpha}\lambda^2}{2}+\frac{\lambda^2}{\nu^{\alpha}})\|u\|^2_{L_2}+2\int_{\R^2}(u_t-\lambda
u)({\mathcal L}u-\lambda u)\le \frac{2}{\mu^{\alpha}}\int_{\R^2}f^2.
\end{equation}
To estimate the second term on the left-hand side of (\ref{eq10}), we use again
the Plansherel's identity to obtain
$$
\int_{\R^2}(u_t-\lambda u)({\mathcal L}u-\lambda
u)=Re[\int_{\R^2}(\lambda+|\omega|^{\alpha})(\lambda+i\tau)|F(u)|^2]=
$$
$$
\int_{\R^2}(\lambda+|\omega|^{\alpha})(\lambda)|F(u)|^2\ge
\int_{\R^2}\lambda^2|F(u)|^2=\lambda^2\|u\|^2_{L_2}\ge 0.
$$ 
We have shown that
$$
(\frac{\mu^{\alpha}\lambda^2}{2}+\frac{\lambda^2}{\nu^{\alpha}}+\lambda^2)\|u\|^2_{L_2}\le
\frac{2}{\mu^{\alpha}}\|f\|^2_{L_2},
$$
or
\begin{equation}\label{eq11}
\|u\|_{L_2}\le M\|f\|_{L_2},
\end{equation}
where the constant $M$ depends on $\mu,\nu,K$, and $\alpha$. 

Obviously,
$$
\|{\mathcal L}u\|_{L_2}\le \|{\mathcal L}u-\lambda u\|_{L_2}+\lambda\|u\|_{L_2},
$$
and
$$
\|u_t\|_{L_2}\le \|u_t-\lambda u\|_{L_2}+\lambda\|u\|_{L_2}
$$
so that the estimate (\ref{eq4}) follows then from (\ref{eq11}), the inequality (\ref{eq50}) and the
established fact that the second term on the left-hand side of (\ref{eq5}) is
non-negative.
$\Box$

\begin{Cor}
\label{supnorm} Under conditions of Lemma \ref{estimate1}, it holds
\begin{equation}
\label{eq12}
\sup_{t,x}|u(t,x)|\le M\|f\|_{L_2}.
\end{equation}
\end{Cor}
\Proof Using the Fourier inversion formula and the Cauchy-Schwarz inequality, we estimate
$$
|u(t,x)|^2\le\Bigl(\frac{1}{(2\pi)^2}\int_{\R^2}|F(u)|d\tau d\omega\Bigr)^2
$$
$$
=\frac{1}{16\pi^4}\Bigl(\int_{\R^2}|F(u)|\Bigl(|-2\lambda+i\tau-|\omega|^{\alpha}|\Bigr)\Bigl(|-2\lambda+i\tau-|\omega|^{\alpha}|\Bigr)^{-1}d\tau d\omega\Bigr)^2\le 
$$
$$
\frac{1}{16\pi^4}I_1I_2,
$$
where
$$
I_1=\int_{\R^2}|F(u)|^2|-2\lambda+i\tau-|\omega|^{\alpha}|^2d\tau d\omega
$$
and
$$
I_2=\int_{\R^2}|-2\lambda+i\tau-|\omega|^{\alpha}|^{-2}d\tau d\omega.
$$
Since $\alpha\in(1,2)$, it follows that
$$
I_2=\int_{\R^2}\frac{d\tau d\omega}{\tau^2+(2\lambda+|\omega|^{\alpha})^2}=\pi\int_{\R}\frac{d\omega}{2\lambda+|\omega|^{\alpha}}:=M_1<\infty.
$$
The term $I_1$ can be estimated as
$$
I_1\le 2\int_{\R^2}|F(u)|^2|-\lambda+i\tau|^2d\tau d\omega+2\int_{\R^2}|F(u)|^2|-\lambda-|\omega|^{\alpha}|^2d\tau d\omega=
$$
$$
2\int_{\R^2}|F(u_t-\lambda u)|^2d\tau d\omega+2\int_{\R^2}|F({\mathcal L}u-\lambda u)|^2d\tau d\omega
$$
$$
=2\|u_t-\lambda u\|^2_{L_2}+2\|{\mathcal L}u-\lambda u\|^2_{L_2}.
$$
Thus, we have shown that
$$
|u(t,x)|^2\le \frac{M_1}{8\pi^4}\Bigl(\|u_t-\lambda u\|^2_{L_2}+\|{\mathcal L}u-\lambda u\|^2_{L_2}\Bigr)
$$
for all $(t,x)\in\R^2$.  The estimate (\ref{eq12}) then follows because of (\ref{eq4}). $\Box$

\section{ Some integral estimates}
\setcounter{equation}{0}

\quad Now, using the analytic estimates from the previous section, we are going to derive the corresponding  integral estimates of Krylov type for the solutions $X$ of the stochastic equation (\ref{eq1}). 

\quad Assume $f\in C_0^{\infty}(\R^2)$ and the coefficients $a$ and $b$ satisfy the assumption (\ref{eq2}). It follows then (see the Appendix)  that the equation (\ref{eq3}) has a solution $u\in H(\R^2)$.  

\quad Let $\psi(t,x)\in C_0^{\infty}(\R^2)$ be a non-negative function with $\psi(t,x)=0$ for all $(t,x)$ such that $|t|+|x|\ge 1$ and $\int_{\R^2}\psi(t,x)dtdx=1$. For $\varepsilon>0$, we define
\[
\psi^{(\varepsilon)}(t,x)=\frac{1}{\varepsilon^2}\psi\Bigl(\frac{t}{\varepsilon},\frac{x}{\varepsilon}\Bigr)
\]
and let $u^{(\varepsilon)}$ to be the convolution of $u$ with the smooth kernel $\psi^{(\varepsilon)}$:
\[
u^{(\varepsilon)}(t,x)=\int_{\R^2}u(s,y)\psi^{(\varepsilon)}(t-s,x-y)dsdy.
\]

\quad Clearly, $u^{(\varepsilon)}\in C^{\infty}_0(\R^2)$ and $\int_{\R^2}\psi^{(\varepsilon)}(s,x)dsdx=1$. Moreover,  $u^{(\varepsilon)}\to u$ as $\varepsilon\to 0$ pointwise and in $L_2(\R^2)$. We also define 
\[
u^{(\varepsilon)}_t:=\frac{\partial}{\partial t}\Bigl(u^{(\varepsilon)}\Bigr) \mbox { and }  u^{(\varepsilon)}_x:=\frac{\partial}{\partial x}\Bigl(u^{(\varepsilon)}\Bigr).
\]

Now, for $\varepsilon>0$, let
\[
f^{(\varepsilon)}:=u^{(\varepsilon)}_t+|b|^{\alpha}{\cal L}u^{(\varepsilon)}+au^{(\varepsilon)}_x-\lambda(1+|b|^{\alpha})u^{(\varepsilon)}.
\]

Because of (\ref{eq3}), $f^{(\varepsilon)}\to f$ as $\varepsilon\to 0$ pointwise and in $L_2(\R^2)$.

\begin{Th}\label{Krylov1} Let $X$ be a solution of the equation (\ref{eq1}) and $\alpha\in(1,2)$. Then, for $t\ge 0,x\in\R$, any measurable function $f:[0,\infty)\times\R\to[0,\infty)$, and $\lambda\ge \delta$, it holds
\begin{equation}\label{eq14}
{\bf E}\int_0^{\infty}f(t,x+X_s)ds\le M\|f\|_{L_2}
\end{equation}	
where the constant $M$ depends on $\nu,\mu,K,$ and $\alpha$.	
\end{Th}	

\Proof  Let $\phi_t=\int_0^t(1+|b(s,X_s)|^{\alpha})ds$.  Then, for all $(t,x)\in[0,\infty)\times\R$, we apply It{\'o}'s formula to the function $u^{(\varepsilon)}(t,X_t)e^{-\phi_t}$ to obtain
\[
{\bf E}u^{(\varepsilon)}(t,X_t)e^{-\lambda\phi_t}-u^{(\varepsilon)}(0,x)=
\]\textbf{}
\[
{\bf E}\int_0^te^{-\lambda\phi_s}\Bigl\{u_t^{(\varepsilon)}(s,X_s)+|b(s,X_s)|^{\alpha}{\mathcal L}u^{(\varepsilon)}(s,X_s)+a(s,X_s)u_x^{(\varepsilon)}(s,X_s)-
\]
\[
\lambda(1+|b|^{\alpha}(s,X_s))u^{(\varepsilon)}(s,X_s)\Bigr\}ds=
-{\bf E}\int_0^te^{-\lambda\phi_s}f^{(\varepsilon)}(s,X_s)ds.
\]
The Corollary \ref{supnorm} implies that $u$ is a bounded function so that the sequence of functions $u^{(\varepsilon)}, \varepsilon>0$ is uniformly bounded. It also follows from its definition and the condition (\ref{eq2}) that the sequence of functions $f^{(\varepsilon)}, \varepsilon>0$ is then also uniformly bounded.

\quad Using the Lebesgue dominated convergence theorem and letting $\varepsilon\to 0$ in the relation
\[
{\bf E}\int_0^te^{-\lambda\phi_s}f^{(\varepsilon)}(s,X_s)ds=u^{(\varepsilon)}(0,x)-{\bf E}u^{(\varepsilon)}(t,X_t)e^{-\lambda\phi_t},
\]
we obtain that
\[
{\bf E}\int_0^te^{-\lambda\phi_s}f(s,X_s)ds=u(0,x)-{\bf E}u(t,X_t)e^{-\lambda\phi_t}.
\]

The above implies
\[
{\bf E}\int_0^tf(s,X_s)ds\le \sup_{t,x}|u(t,x)|\le M\|f\|_{L_2}.
\]
Using the Fatou's lemma and letting $t\to\infty$, we obtain
\[
{\bf E}\int_0^{\infty}f(s,X_s)ds\le M\|f\|_{L_2}.
\]
The later inequality can be extended to any nonnegative measurable function $f$ by using the standard  arguments of a monotone class theorem (see, for example, Theorem 21 in \cite{DM}). $\Box$

\quad We can also obtain a local version of the estimate (\ref{eq14}). For that, for any $t>0$ and $m\in\N$, we define $\|f\|_{2,m,t}:=(\int_0^t\int_{[-m,m]}|f(s,x)|dsdx)^{1/2}$  as the $L_2$-norm of $f$ on $[0,t]\times[-m,m]$. Let also $\tau_m(X)=\inf\{t\ge 0: |X_t|>m\}$. Then, applying  (\ref{eq14}) to the function $\bar f(s,x)=f(s,x){\bf 1}_{[0,t]\times[-m,m]}(s,x)$, we obtain
\begin{Cor}\label{local1} Let $X$ be a solution of equation (\ref{eq1}) with $\alpha\in(1,2)$ and the assumption (\ref{eq2}) is satisfied. Then, for any $t>0,m\in\N$, and any nonnegative  measurable function $f$, it holds that
\begin{equation}\label{localestimate}	
{\bf E}\int_0^{t\land\tau_m(X)}f(s,X_s)ds\le M\|f\|_{2,m,t}
\end{equation}
where the constant $M$ depends on $\mu,\nu,K,t,$ and $m$.
\end{Cor}

\section{Existence of solutions for stochastic equations with measurable coefficients}
\setcounter{equation}{0}
\quad As an applications of the integral estimates derived in the previous section, we prove here the existence of solutions for the SDE (\ref{eq1}) under the assumption (\ref{eq2}) where $Z$ is a symmetric stable process of index $\alpha\in(1,2]$.

\quad For $\alpha=2$, the existence of solutions under (\ref{eq2}) is well-known (cf. \cite{Krylov}). Hencefore, we restrict ourself to the case $1<\alpha<2$. 
\begin{Th}\label{boundedcoefficients} Assume that $a(t,x)$ and $b(t,x)$ are two measurable functions satisfying the condition (\ref{eq2}) and $\alpha\in(1,2)$. Then, for any $x_0\in \R$, there exists a solution of the equation (\ref{eq1}).
\end{Th}
\Proof Because of the assumptions (\ref{eq2}), for $n=1,2,\dots$, there are sequences of functions $a_n(t,x)$ and $b_n(t,x)$ such that they are globally Lipshitz continuous, uniformly bounded and  $a_n\to a, b_n\to b$ (a.s.) as $n\to\infty$. For any $n=1,2,\dots$, the equation (\ref{eq1}) has a unique solution, even so-called strong solution (see, for example, Theorem 9.1 in \cite{IkWa}). That is, for any fixed symmetric stable process $Z$ defined on a probability space $(\Omega, \mathcal{F}, \PP)$, there exists a sequence of processes $X^n, n=1,2\dots,$ such that
\begin{equation}\label{Eq1}
dX_t^n=b_n(t,X^n_{t-})dZ_t+a_n(t,X_t^n)dt,\quad X_0^n=x_0\in\R,\quad t\ge 0.
\end{equation}
Let
\[
M^n_t:= \int_0^tb_n(s,X^n_{s-})dZ_s \mbox { and } Y^n_t:=\int_0^ta_n(s,X^n_s)ds
\]
so that
\[
X^n=x_0+M^n+Y^n,\quad n\ge 1.
\]

\quad As next, we show that the sequence of processes
$H^n:=(X^n,M^n,Y^n, Z), n\ge 1$, is tight in the sense of weak convergence in $(\D^4,\mathcal{D}^4)$. Due to the well-known Aldous' criterion (\cite{Ald}), it suffices to show that 
\begin{equation}\label{Aldous1}
\lim_{l\to\infty}\limsup_{n\to\infty}\PP\Bigl[\sup_{0\le s\le t}\|H^n_s\|>l\Bigr]=0
\end{equation}
for all $t\ge 0$ and
\begin{equation}\label{Aldous2}
\limsup_{n\to\infty}\PP\Bigl[\|H^n_{t\land(\tau^n+\delta_n)}-H^n_{t\land\tau^n}\|>\varepsilon\Bigr]=0
\end{equation}
for all $t\ge 0, \varepsilon>0$, every sequence of $\F$-stopping times $\tau^n$, and every sequence of real numbers $\delta_n$ such that $\delta_n\downarrow 0$. Here $\|\cdot\|$ denotes the Euclidean norm of a vector.

\quad It is clear that for this it suffices only to verify that the sequence of processes $(M^n,Y^n)$ is tight in $(\D^2,\mathcal{ D}^2)$. But this is trivially fulfilled because of the uniform boundness of the coefficients $a_n$ and $b_n$ for all $n\ge 1$. 

\quad From the tightness of the sequence $\{H^n \}$ we conclude that there exists a subsequence $\{n_k\},k=1,2,\dots$, a probability space $(\bar\Omega,\bar{\mathcal{F}},\bar\PP)$ and the process $\bar H$ on it with values in $(\D^4,{\mathcal{D}}^4)$ such that $H^{n_k}$ converges weakly (in distribution) to the process $\bar H$ as $k\to\infty$. For simplicity, let $\{n_k\}=\{n\}$.  

\quad We use now the well-known embedding principle of Skorokhod (see, e.g. Theorem 2.7 in \cite{IkWa}) to imply the convergence of the sequence $\{H^n\}$ a.s. in the following sense: there exists a probability space $(\tilde\Omega,\tilde{\mathcal{F}},\tilde\PP)$ and processes $\tilde H=(\tilde X,\tilde M,\tilde Y,\tilde Z),\quad \tilde H^n=(\tilde X^n,\tilde M^n,\tilde Y^n,\tilde Z^n),\quad n=1,2,\dots,$ on it such that 
\begin{enumerate}
	\item[1)] $\tilde H^n\to\tilde H$ as $n\to\infty$ $\tilde\PP$-a.s.
	\item[2)] $\tilde H^n=H^n$ in distribution for all $n=1,2,\dots.$
\end{enumerate}
Using standard measurability arguments (\cite{Krylov}, chapter 2), one can prove that the processes $\tilde Z^n$ and $\tilde Z$ are symmetric stable processes of the index $\alpha$ with respect to the augmented filtrations $\tilde\F^n$ and $\tilde\F$ generated by processes $\tilde H^n$ and $\tilde H$, respectively. 

\quad Relying on the above properties $1)$ and $2)$, and the equation (\ref{Eq1}), one can show (\cite{Krylov}, chapter 2) that
\[
\tilde X^n_t=x_0+\int_0^tb_n(s,X^n_{s-})\tilde Z^n_s+\int_0^ta_n(s,\tilde X^n_s)ds, \quad t\ge 0,\quad\tilde\PP\mbox{-a.s.}
\]

At the same time, from the properties $1), 2)$ and the quasi-left continuity of the the processes $\tilde X^n$ it follows that
\begin{equation}\label{convergence1}
\lim_{n\to\infty}\tilde X^n_t=\tilde X_t, \quad t\ge 0,\quad \tilde\PP\mbox{-a.s.}
\end{equation}
Hence in order to show that the process $\tilde X$ is a solution of the equation (\ref{eq1}), it is enough to prove that, for all $t\ge 0$,
\begin{equation}\label{convergence3}
\lim_{n\to\infty}\int_0^tb_n(s,\tilde X^n_s)d\tilde Z^n_s=\int_0^tb(\tilde X_s)d\tilde Zs\quad \tilde\PP\mbox{- a.s. }
\end{equation}
and
\begin{equation}\label{convergence4}
\lim_{n\to\infty}\int_0^ta_n(s,\tilde X^n_s)ds=\int_0^ta(\tilde X_s)ds\quad \tilde\PP\mbox{- a.s. }
\end{equation}
Now we remark that from the convergence in probability it follows that there is a subsequence for which the convergence with probability one holds. Therefore,  to verify (\ref{convergence3}) and (\ref{convergence4}) , it suffices to show that for all $t\ge 0$ and $\varepsilon>0$ we have 
\begin{equation}\label{convergence3'}
\lim_{n\to\infty}\tilde\PP\Bigl[|\int_0^{t}b_n(s,\tilde X^n_s)d\tilde Z^n_s-\int_0^{t}b(s,\tilde X_s)d\tilde Z_s|>\varepsilon\Bigr]=0
\end{equation}
and
\begin{equation}\label{convergence4'}
\lim_{n\to\infty}\tilde\PP\Bigl[|\int_0^{t}a_n(s,\tilde X^n_s)ds-\int_0^{t}a(s,\tilde X_s)ds|>\varepsilon\Bigr]=0.
\end{equation}

We will also need the following result that can be proven in the same way as Lemma 4.2 in \cite{Ku1}.

\begin{Le}\label{limitprocess} Let $\tilde X$ be the process as defined above. Then, for any Borel measurable function $f:[0,\infty)\times\R\to[0,\infty)$ and any $t\ge 0,$ there exists a sequence $m_k\in (0,\infty), k=1,2,\dots$ such that $m_k\uparrow\infty$ as $k\to\infty$ and it holds
	\[
	{\bf \tilde E}\int_0^{t\land\tau_{m_k}(\tilde X)}f(s,\tilde X_s)ds\le M \|f\|_{2,m_k,t},
	\]
	where the constant $M$ depends on $\lambda,\alpha, t$ and $m_k$ only. Moreover, it holds
	\begin{equation}\label{help}
	\tilde \PP\Bigl[\tau_m(\tilde X^n)<t\Bigr]\to\tilde\PP\Bigl[\tau_m(\tilde X)<t\Bigr]\mbox { as } n\to\infty.
	\end{equation}
\end{Le}

\quad Without loss of generality, we can assume $\{m_k\}=\{m\}$. 

\quad Let us prove (\ref{convergence3'}) and (\ref{convergence4'}). For a fixed $k_1\in\N$ we have
\[
\tilde \PP\Bigl[|\int_0^{t}b_n(s,\tilde X^n_{s-})d\tilde Z^n_s-\int_0^{t}b(s,\tilde X_{s-})d\tilde Z_s|>\varepsilon\Bigr]\le
\]
\[
\tilde \PP\Bigl[|\int_0^{t}b_{k_1}(s,\tilde X^n_{s-})d\tilde Z^n_s-\int_0^{t}b_{k_1}(s,\tilde X_{s-})d\tilde Z_s|>{\varepsilon\over 3}\Bigr]
\]
\[
+ \tilde \PP\Bigl[|\int_0^{t\land\tau_m(\tilde X^n)}b_{k_1}(s,\tilde X^n_s)d\tilde Z^n_s-\int_0^{t\land\tau_m(\tilde X^n)}b_n(s,\tilde X^n_{s-})d\tilde Z^n_s|>{\varepsilon\over 3}\Bigr]
\]
\[
+ \tilde \PP\Bigl[|\int_0^{t\land\tau_m(\tilde X)}b_{k_1}(s,X_s)d\tilde Z_s-\int_0^{t\land\tau_m(\tilde X)}b(s,\tilde X_{s-})d\tilde Z_s|>{\varepsilon\over 3}\Bigr]
\]
\[
+\tilde\PP\Bigl[\tau_m(\tilde X^n)<t\Bigr]+\tilde\PP\Bigl[\tau_m(\tilde X)<t\Bigr].
\]
The first term on the right side of the inequality above converges to $0$ as $n\to\infty$ by Chebyshev's inequality and Skorokhod lemma for stable integrals (see \cite{PragZa}, Lemma 2.3). To show the convergence to $0$ as $n\to\infty$ of the second and third terms we use first the Chebyshev's inequality and then Corollary \ref{localestimate} and Lemma \ref{limitprocess}, respectively. We obtain
\[
\tilde \PP\Bigl[|\int_0^{t\land\tau_m(\tilde X^n)}b_{k_1}(s,\tilde X^n_s)d\tilde Z^n_s-\int_0^{t\land\tau_m(\tilde X^n)}b_n(s,\tilde X^n_{s-})d\tilde Z^n_s|>{\varepsilon\over 3}\Bigr]
\]
\begin{equation}\label{est4}
\le \frac{3}{\epsilon}\tilde{\bf E}|\int_0^{t\land\tau_m(s,\tilde X^n)}|b_{k_1}-b_n|^{\alpha}(s,\tilde X^n_{s-})ds|\le {3\over \varepsilon}M\||b_{k_1}-b_n|^{\alpha}\|_{2,m,t}
\end{equation}
and
\[
\tilde \PP\Bigl[|\int_0^{t\land\tau_m(\tilde X)}b_{k_1}(s,X_s)d\tilde Z_s-\int_0^{t\land\tau_m(\tilde X)}b(s,\tilde X_{s-})d\tilde Z_s|>{\varepsilon\over 3}\Bigr]
\]
\begin{equation}\label{est5}
\le\frac{3}{\epsilon}\tilde {\bf E}|\int_0^{t\land\tau_m(\tilde X)}|b_{k_1}-b|^{\alpha}(s,\tilde X_{s-})ds|\le {3\over \varepsilon}M\||b_{k_1}-b|^{\alpha}\|_{2,m,t}
\end{equation}
where the constant $M$ depends on $\mu,\nu, K,m,t,$ and $\alpha$ only. 

It follows from the definition of the sequence $b_n$ that, for any $m\in\N$, $|b_{k_1}-b|^{\alpha}\to 0$ as $k_1\to\infty$ in $L_{2,m,t}$-norm. Then, passing to the limit in (\ref{est4}) and (\ref{est5}) first $n\to\infty$ and then $k_1\to\infty$, we obtain that the right sides of (\ref{est4}) and (\ref{est5}) converge to $0$. 

\quad Because of the property (\ref{help}), 
the remaining terms $\tilde \PP\Bigl[\tau_m(\tilde X^n)<t\Bigr]$ and  $\tilde\PP\Bigl[\tau_m(\tilde X)<t\Bigr]$ can be made arbitrarily small by choosing large enough $m$ for all $n$ due to the fact that the sequence of processes $\tilde X^n$ satisfies the property (\ref{Aldous1}).
This verifies (\ref{convergence3'}). The convergence (\ref{convergence4'}) can be verified similarly. We omit the details. 

\quad Thus, we have proven the existence of the process $\tilde X$ that solves the equation (\ref{eq1}).  $\Box$
	
\section{Appendix} 
\setcounter{equation}{0}

\quad Here we prove the existence of a solution of equation (\ref{eq3}) in the Sobolev space $H(\R^2)$ for any $f\in L_2(\R^2)$ and the coefficients $a$ and $b$ satisfying the condition (\ref{eq2}). In order to do so, we use the method of continuity and the method of a priori estimates in a similar way as it is done in \cite{Krylov1} in the case of classical elliptic and  parabolic equations. 

\quad We first start with the equation
\begin{equation}\label{simpleeq}
u_t+{\cal L}u-\lambda u=f,
\end{equation}
where $\lambda>0$ and $f\in L_2(\R^2)$ is a given function.

\quad To solve (\ref{simpleeq})  in $H(\R^2)$, we will need some lemmas.
\begin{Le}\label{lemma1} Let $u\in C_0^2(\R^2)$ be a solution of (\ref{simpleeq}). Then, it holds
\begin{equation}\label{eq20}
\|u_t\|_{L_2}^2+\lambda^2\|u\|_{L_2}^2+\|{\cal L}u\|_{L_2}^2\le \|f\|_{L_2}^2.
\end{equation}		
\end{Le}	
\Proof Applying the Fourier transform in variables $(t,x)$ to the equation (\ref{simpleeq}), we obtain
\[
-i\tau F[u]-(\lambda+|w|^{\alpha})F[u]=F[f],
\]
or,
\[
\Bigl(|\tau|^2+(\lambda+|w|^{\alpha})^2\Bigr)|F[u]|^2=|F[f]|^2
\]
which implies
\[
|\tau|^2|F[u]|^2+\lambda^2|F[u]|^2+|w|^{2\alpha})|F[u]|^2\le |F[f]|^2.
\]

Integrating the last relation over $\R^2$ and using the Parseval's identity, we obtain (\ref{eq20}). $\Box$

\begin{Le}\label{lemma2} Let $\lambda>0$ and $u\in C^{\infty}_0(\R^2)$ so that
\[
u_t+{\cal L}u-\lambda u=0.
\]
Then $u=0$ a.e.
\end{Le}
\Proof As in Lemma \ref{lemma1}, we apply the Fourier transform to (\ref{simpleeq}) (with $f=0$), to get
\[
\|u_t\|_{L_2}^2+\lambda^2\|u\|_{L_2}^2+\|{\cal L}u\|_{L_2}^2\le 0.
\]
It follows then that $\|u\|_{L_2}=0$ and $u=0$ a.e. in $\R^2$.$\Box$

\quad Let $\lambda>0$ and consider a set of functions
\[
{\cal A}:=\{g(t,x)=\frac{\partial}{\partial t}h(t,x)+{\cal L}h(t,x)-\lambda h(t,x)\mbox { for some } h\in C_0^{\infty}(\R^2)\}.
\]

\begin{Le}\label{lemma3} The set ${\cal A}$ is dense in $L_2(\R^2)$.	
\end{Le}	
	
\Proof From the converse. If ${\cal A}$ is not dense in $L_2(\R^2)$, then by the Hahn-Banach theorem there is a function $g\in L_2(\R^2)$  with $\|g\|_{L_2}\neq 0$ so that
\[
\int_{R^2}g(t,x)\Bigl(\frac{\partial}{\partial t}+{\cal L}-\lambda\Bigr)u(t,x)dtdx=0
\]
for all $u\in C_0^{\infty}(\R^2)$.

The last relation also implies that
\begin{equation}\label{eq21}
\int_{R^2}g(t,x)\Bigl(\frac{\partial}{\partial t}+{\cal L}-\lambda\Bigr)u(\tau-t,y-x)dtdx=0
\end{equation}
since $u(\tau-t,y-x)\in C_0^{\infty}(\R^2)$ for all fixed $(\tau,y)\in\R^2$ .

\quad Using convolution, (\ref{eq21}) is then written as
\begin{equation}\label{eq22}
g\star \frac{\partial }{\partial t}u(\tau, y)+g\star{\cal L}u(\tau,y)-\lambda g\star u(\tau,y)=0.
\end{equation}

Clearly, 
\begin{equation}\label{eq23}
g\star \frac{\partial }{\partial t}u=\frac{\partial }{\partial t}\Bigl(g\star u\Bigr).
\end{equation}
We also have that
\[
g\star{\cal L}u(\tau,y)=\int_{\R^2}g(t,x){\cal L}u(\tau-t,y-x)dtdx=
\]
\[
\int_{\R^2}g(t,x)\int_{\R}\Bigl[u(\tau-t,y-x+z)-u(\tau-t,y-x)-{\bf 1}_{|z|<1}u_x(\tau-t,y-x)z\Bigr]\frac{dz}{|z|^{1+\alpha}}dtdx
\]
and
\[
{\cal L}(g\star u)(\tau,y)=\int_{\R}\Bigl(\int_{\R^2}g(t,x)u(\tau-t,y-x+z)dtdx-
\]
\[
-\int_{\R^2}g(t,x)u(\tau-t,y-x)dtdx-\int_{\R^2}zg(t,x)u_x(\tau-t,y-x){\bf 1}_{|z|<1}dtdx\Bigr)\frac{dz}{|z|^{1+\alpha}}=
\]
\[
\int_{\R^2}g(t,x)\int_{\R}\Bigl[u(\tau-t,y-x+z)-u(\tau-t,y-x)-{\bf 1}_{|z|<1}u_x(\tau-t,y-x)z\Bigr]\frac{dz}{|z|^{1+\alpha}}dtdx
\]
where we used the fact that $(g\star u)_x=g\star u_x$.

\quad Comparing the above relations we conclude that
\begin{equation}\label{eq24}
g\star{\cal L}u={\cal L}(g\star u).
\end{equation}
Using (\ref{eq23}) and (\ref{eq24}), the equation (\ref{eq22}) becomes
\[
\Bigl(\frac{\partial}{\partial t}+{\cal L}-\lambda\Bigr)g\star u(\tau,y)=0.
\]
We then apply Lemma \ref{lemma2} to obtain
\[
\int_{\R^2}g(t,x)u(\tau-t,y-x)dtdx=0
\]
for all $u\in C_0^{\infty}(\R^2)$ and all $(\tau,y)\in\R^2$. It follows from the general integration theory that $g=0$ a.e. in $\R^2$ implying $\|g\|_{L_2}=0$ which is a contradiction. $\Box$

\begin{Le}\label{lemma4} Let $\lambda>0$ and $f\in L_2(\R^2)$. Then, there is a solution $u\in H(\R^2)$  of the equation (\ref{simpleeq}).
\end{Le}	
\Proof By Lemma \ref{lemma3}, there is a sequence of functions $u^n\in C_0^{\infty}(\R^2)$ so that
\[
\Bigl(u^n_t+{\cal L}u^n-\lambda u^n\Bigr)\to f\mbox{ as } n\to\infty
\]
in $L_2(\R^2)$.

\quad Define
\begin{equation}\label{eq25}
f^n:=\Bigl(u^n_t+{\cal L}u^n-\lambda u^n\Bigr), n=1,2,...
\end{equation}

\quad Using Lemma \ref{lemma1}, we obtain that
\[
\|u_t^n-u_t^m\|_{L_2}^2+\lambda^2\|u^n-u^m\|_{L_2}^2+\|{\cal L}u^n-{\cal L}u^m\|_{L_2}^2\le \|f^n-f^m\|_{L_2}^2
\]
for all $n,m=1,2,...$

\quad Since $(f^n)$ converges in $L_2(\R^2)$, it is a Cauchy sequence so that $\|f^n-f^m\|_{L_2}\to 0$ as $n,m\to\infty$. This implies that the sequences $(u^n), (u^n_t)$, and $({\cal L}u^n)$ are also Cauchy sequences. Because of the completeness of $L_2(\R^2)$, the following limits will exist in $L_2(\R^2)$:
\[
v(t,x):=\lim_{n\to\infty}u^n(t,x), v_t(t,x):=\lim_{n\to\infty}u^n_t(t,x), {\cal L}v(t,x):=\lim_{n\to\infty}{\cal L} u^n(t,x).
\] 

It follows then from (\ref{eq25}) that it holds
\[
v_t+{\cal L} v-\lambda v=f \mbox { a. e.  in } \R^2.
\]
Therefore, $v$ is a solution of the equation (\ref{simpleeq}) in the sense described above which is often referred to as {\it a generalized solution in the Sobolev space $H$}.
 $\Box$

\quad Now, for any $\lambda>0$ and $\alpha\in(1,2)$, we consider the operator
\[
L:=\frac{\partial}{\partial t}+|b|^{\alpha}{\mathcal L} +a\frac{\partial}{\partial x}-\lambda(1+|b|^{\alpha}),
\]
where the real-valued functions $a(t,x), b(t,x)$ satisfy the assumption (\ref{eq2}).

\quad For any $s\in[0,1]$, we set
\[
L_s:=(1-s)(\frac{\partial}{\partial t}+{\mathcal L} -\lambda)+sL.
\]

The following result is the analog of Theorem 1.4 from \cite{Krylov1}. The proof is entirely based on general functional analysis facts and we refer for details to \cite{Krylov1}. 
\begin{Prop}\label{continuity} Assume that there are constants $\lambda>0$ and $M\in(0,\infty)$ such that for any $u\in C^2_0(\R^2)$ and $s\in[0,1]$ it holds
\begin{equation}\label{eq30}
\|u\|_H\le M\|L_s\|_{L_2}.	
\end{equation}
Then, for any $f\in L_2(\R^2)$, there is a function $u\in H(\R^2)$ satisfying	$L u=f$.
\end{Prop}	
The condition (\ref{eq30}) can be reformulated in the following form: for any $u\in H(\R^2)$ satisfying the equation $L_s u=f$, it holds
\begin{equation}\label{eq31}
\|u\|_H\le M\|f\|_{L_2}.
\end{equation} 

The estimate  (\ref{eq31}) is called {\it an a priori estimate} for the equation $L_{s}u=f$ since we do not know the existence of such a solution yet.

\begin{Prop}\label{lemma5} For any $f\in L_2(\R^2)$, there is a solution $u\in H(\R^2)$ of the equation $Lu=f$.
	
\end{Prop}	

\Proof Step 1. Assume first that $a=0$. It follows from Lemma \ref{estimate1} that, for any $u\in C_0^2(\R^2)$ and any $\lambda>0$, it holds 
\begin{equation}\label{eq35}
\|u_t\|_{L_2}^2+\lambda^2\|u\|_{L_2}^2+\|{\cal L}u\|_{L_2}^2\le M\|u_t+|b|^{\alpha}{\mathcal L}u -\lambda(1+|b|^{\alpha})u\|_{L_2}^2,
\end{equation}
where the constant $M$ depends on $\nu$ and $\mu$.

\quad For $s\in[0,1]$, we consider
\[
\tilde L_su:=(1-s)(u_t+{\mathcal L}u -2\lambda u)+s\Bigl( u_t+|b|^{\alpha}{\mathcal L}u-\lambda(1+|b|^{\alpha})u\Bigr).
\]
It can easily be seen that
\[
\tilde L_s=u_t+[1-s+s|b|^{\alpha}]{\mathcal L}u-\lambda[1+1-s+s|b|^{\alpha}]u=
\]
\[
u_t+\sigma(s){\mathcal L}u-\lambda[1+\sigma(s)]u
\]
where
\[
\sigma(s)=1-s+s|b|^{\alpha}.
\]
\quad Because of Lemma \ref{lemma4} , the equation $u_t+{\mathcal L}u -2\lambda u=f$ has a solution $u\in H(\R^2)$ for any $\lambda>0$ and $f\in L_2(\R^2)$. Using Proposition \ref{continuity}, the assertion in Step 1 is then proved if, for any $s\in[0,1]$ and any $u\in C^2_0(\R^2)$, it follows that
\[
\|u\|_H\le M\|\tilde L_s u\|_{L_2}.
\] 
The later, however, follows from (\ref{eq35}) if we replace $|b|^{\alpha}$ by $\sigma(s)$ and notice that, for any $s\in[0,1]$, it holds 
\[
0<\min\{1,\mu^{\alpha}\}\le \sigma(s)\le \max\{1,\nu^{\alpha}\}
\]
since $\sigma(s)$ is a linear function in $s$.

\quad Step 2. For $s\in[0,1]$, we consider the operator
\[
L_s=(1-s)\Bigl( u_t+|b|^{\alpha}{\mathcal L}u-\lambda(1+|b|^{\alpha})u\Bigr)+s Lu=
\]
\[
u_t+|b|^{\alpha}{\mathcal L}u-\lambda(1+|b|^{\alpha})u+s au_x.
\]
Using (\ref{eq35}), we obtain that, for any $u\in C^2_0(\R^2)$ and $\lambda>0$
\begin{equation}\label{eq40}
\|u_t\|_{L_2}+\lambda\|u\|_{L_2}+\|{\cal L}u\|_{L_2}\le M_1\|L_s\|_{L_2}+M_2\|u_x\|_{L_2}
\end{equation}
where the constants $M_1$ and $M_2$ depend on the bounds of the coefficients $a$ and $b$.

\quad It can be easily seen that, for any fixed  $1<\alpha<2$, there exists $\lambda_0>0$ so that
\[
M_2|\omega|^2\le \frac{1}{2}(\lambda_0+|\omega|^{\alpha})^2, \omega\in\R.
\]
\quad It follows then that
\[
M_2\|u_x\|_{L_2}\le \frac{1}{2}\|{\cal L}u\|_{L_2}+\frac{\lambda_0}{2}\|u\|_{L_2}
\]
and using (\ref{eq40}) we conclude that
\[
\|u_t\|_{L_2}+(\lambda-\frac{\lambda_0}{2})\|u\|_{L_2}+\frac{1}{2}\|{\cal L} u\|_{L_2}\le M_1\|L_s u\|_{L_2}.
\]
The last relation implies a priori estimate
\[
\|u\|_H\le M\|L_s u\|_{L_2}
\]
for $\lambda>\lambda_0/2$ with $M$ depending on the bounds of $a$ and $b$. The later, in turn, implies the existence of a solution $u\in H(\R^2)$ of the equation $Lu=f$ for any $f\in L_2(\R^2)$ because of Proposition \ref{continuity}. $\Box$

\end{document}